\documentclass[12pt]{amsart}

\usepackage{longtable}
\usepackage{url}
\usepackage{enumerate}
\usepackage{epsfig}

 \usepackage{ifthen}
 \newcommand{\mymarginpar}[1]{%
    \marginpar{\ifthenelse{\isodd{\arabic{page}}}{\flushleft #1}{\flushright #1}}}

 \renewcommand{\phi}{\varphi}
 \newcommand{\alg}[1]{\mathcal{#1}}                % algebras
 \newcommand{\abs}[1]{\left\lvert #1 \right\rvert}
 \newcommand{\norm}[1]{\left\lVert #1 \right\rVert}
 \newcommand{\snorm}[1]{\lVert #1 \rVert}

 \newcommand{\Bignorm}[1]{\Bigl\lVert #1 \Bigr\rVert}

 \DeclareMathOperator{\linspan}{{\mathrm span}}                   % span
 \DeclareMathOperator{\diam}{{\mathrm diam}}                   % support
 \DeclareMathOperator{\supp}{{\mathrm supp}}                   % support
 \renewcommand{\Im}{\mathop{\mathrm{Im}}}
 \newcommand{\IC}{\mathbf{C}}                     % komplexe Zahlen
 \newcommand{\ID}{\mathbf{D}}                     % unit disk
 \newcommand{\IF}{\mathbf{F}}                     % free group
 \newcommand{\IZ}{\mathbf{Z}}                     % ganze Zahlen
 \newcommand{\eps}{{\varepsilon}}                 % epsilon

 \theoremstyle{plain} %%%%%%%%%%%%%%%%%%%%%%%%%%%%%%%%%
 \newtheorem{Theorem}{Theorem}[section]
 \newtheorem{Lemma}[Theorem]{Lemma}
 \newtheorem{Proposition}[Theorem]{Proposition}
 
 \newtheorem{Corollary}[Theorem]{Corollary}
 
 \theoremstyle{definition} %%%%%%%%%%%%%%%%%%%%%%%%%%%%
 \newtheorem{Definition}[Theorem]{Definition}
 
 \newtheorem{Example}[Theorem]{Example}

 \numberwithin{equation}{section}

\begin{document}
 \title[Free spectra]{On the computation of spectra in free probability}
 \author{Franz Lehner}
 \address{
 In\-sti\-tut f\"ur Ana\-ly\-sis und Nu\-me\-rik\\
 Jo\-han\-nes Kep\-ler Uni\-ver\-sit\-\"at Linz\\
 %Al\-ten\-ber\-ger\-str. 79\\
 4040 Linz\\
 Austria
 }
 \curraddr{\'Equipe d'Analyse\\
   Universit\'e Paris 6\\
   Bo\^te 186\\
   Tour 46-0 \\
   4, place Jussieu\\
   75252 Paris Cedex 05
 }
 \email{lehner@bayou.uni-linz.ac.at}
 \date{\today}


 \begin{abstract}
   We use free probability techniques to compute borders of spectra of
   non hermitian operators in finite von Neumann algebras which arise as
   ``free sums'' of ``simple'' operators.
   To this end, the resolvent is analyzed with the aid of the Haagerup inequality.
   Concrete examples coming from reduced $C^*$-algebras of free product
   groups and leading to systems of polynomial equations illustrate the
   approach.
 \end{abstract}

 \maketitle{}

 %SSSSSSSSSSSSSSSSSSSSSSSSSSSSSSSSSSSSSSSSSSSSSSSSSSSSSSSSSSSSSSSSSSSSSS
 %SSSSSSSSSSSSSSSSSSSSSSSSSSSSSSSSSSSSSSSSSSSSSSSSSSSSSSSSSSSSSSSSSSSSSS
 \section{Introduction}
 \label{sec:introduction}
 %SSSSSSSSSSSSSSSSSSSSSSSSSSSSSSSSSSSSSSSSSSSSSSSSSSSSSSSSSSSSSSSSSSSSSS
 %SSSSSSSSSSSSSSSSSSSSSSSSSSSSSSSSSSSSSSSSSSSSSSSSSSSSSSSSSSSSSSSSSSSSSS

 Computation of spectra of convolution operators on discrete groups
  was one of the motivations behind the
 development of free probability.  It has some interest for harmonic
 analysis and analysis of random walks on free
 products of discrete groups, where under the name of the
 \emph{transition operator} these operators
 carry much information about the random walks under consideration.
 \emph{Free probability} is an abstract framework for
 harmonic analysis on the free group in the language of non-commutative
 probability. Let us introduce briefly the terminology of the latter,
 which is used throughout the paper.

 \begin{Definition}
   A \emph{non-commutative probability space} is a pair $(\alg{A},\phi)$ of a
   (complex) algebra $\alg{A}$ with unit $I$ and a linear functional
   $\phi$ on $\alg{A}$ which satisfies $\phi(I)=1$.  We will usually
   work with $C^*$-algebras and faithful states, i.e.\ positive unital linear
   functionals. The elements of the algebra are called
   \emph{(non-commutative) random variables}. The distribution of a
   random variable $a\in\alg{A}$  is given by the collection of its
   \emph{moments} $\phi(a^n)$, $n=0,1,2,\dots$. If $a$ is self-adjoint,
   then this corresponds to a probability measure on the spectrum of $a$.
   The distribution of a family of random variables is the collection
   of their mixed moments, which abstractly can be interpreted as a
   linear functional on the polynomials on noncommutative variables.

   Given a noncommutative probability space $(\alg{A},\phi)$, the
   subalgebras $\alg{A}_i\subseteq \alg{A}$ are called \emph{freely
   independent} (or \emph{free} for short) if
   \begin{equation}
     \label{eq:FreenessDefinition}
     \phi(a_1a_2\cdots a_n)=0
   \end{equation}
   whenever $a_j\in \alg{A}_{i_j}$ with
   $\phi(a_j)=0$ and $i_j\ne i_{j+1}$ for $j=1,\dots,n-1$.
 \end{Definition}

 In the rest of this section we review some of the necessary facts from
 free probability theory and refer the reader to
 \cite{Voiculescu.Dykema.Nica:1992:free} and
 \cite{Hiai.Petz:2000:semicircle} for further information.

 One of the basic problems is \emph{free convolution}:  Given free
 non-commutative random variables $a$, $b\in\alg{A}$ whose moments are
 known, compute the moments of their sum $a+b$.
 From \eqref{eq:FreenessDefinition} one sees
 immediately by induction that the mixed moments of $a$ and $b$ only
 depend on the individual distributions of $a$ and $b$.  Therefore the
 moments of $a+b$ only depend on the individual moments of $a$ and $b$,
 and thus the term \emph{free convolution} of distributions is
 justified.  In the case of selfadjoint random variables, once the
 moments are known, one can proceed to compute the spectra, see
 section~\ref{sec:SelfadjointOperators}.  What to do in the
 non-selfadjoint case will be the subject of the rest of this paper.

 The computational machinery of free convolution was found independently and
 about the same time by W.~Woess \cite{Woess:1986:nearest},
 D.~Cartwright and P.~Soardi \cite{Cartwright.Soardi:1986:random},
 J.C.~McLaughlin \cite{McLaughlin:1986:random} in the language of
 random walks, and in most generality by D.~Voiculescu
 \cite{Voiculescu:1986:addition}. Using the conventions of the latter,
 the recipe goes as follows.  The Cauchy transform of a random variable
 $a$ is the function
 \begin{equation}
   \label{eq:CauchyTransform}
   G_a(\zeta) = \phi( (\zeta-a)^{-1} ) = \frac{1}{\zeta}
   \sum_{n=0}^\infty
    \frac{\phi(a^n)}{\zeta^n}
 \end{equation}
 which is defined and analytic at least for $\abs{\zeta}>\norm{a}$.  It
 has an inverse (under composition) in some neighbourhood of infinity
 which has the form
 \begin{equation}
   \label{eq:Ktransform}
   K_a(z) = G_a^{-1}(z) = \frac{1}{z}(1 + R_a(z))
 \end{equation}
 where $R_a(z)= c_1 z + c_2 z^2+ \dots$ is analytic and one has
 \begin{equation}
   \label{eq:RTransformSum}
   R_{a+b}(z) = R_a(z) + R_b(z)
   .
 \end{equation}
 This allows in principle to compute the moments of the sum $a+b$.
 However, function inversion is difficult and for non-selfadjoint
 operators knowledge of moments is not sufficient to determine the
 spectrum. For the latter, a more detailed analysis of the resolvent is
 necessary.

 The paper is organized as follows.

 In section~\ref{sec:SelfadjointOperators} we review the case of
 self-adjoint operators, where the computation is reduced to the
 solution of a moment problem.

 In section~\ref{sec:nons.a.} we use Haagerup inequality to estimate
 the norm of the resolvent of a sum of free operators.

 In section~\ref{sec:TwoOperatorResolvent} we look at the sum of two
 free operators and give a somewhat easier proof for this case.

 Finally in section~\ref{sec:Examples} we present some examples of
 computations. 

 %SSSSSSSSSSSSSSSSSSSSSSSSSSSSSSSSSSSSSSSSSSSSSSSSSSSSSSSSSSSSSSSSSSSSSS
 %SSSSSSSSSSSSSSSSSSSSSSSSSSSSSSSSSSSSSSSSSSSSSSSSSSSSSSSSSSSSSSSSSSSSSS
 \section{Self-adjoint operators}
 \label{sec:SelfadjointOperators}
 %SSSSSSSSSSSSSSSSSSSSSSSSSSSSSSSSSSSSSSSSSSSSSSSSSSSSSSSSSSSSSSSSSSSSSS
 %SSSSSSSSSSSSSSSSSSSSSSSSSSSSSSSSSSSSSSSSSSSSSSSSSSSSSSSSSSSSSSSSSSSSSS

 We collect here some well known techniques about computation of
 spectra of sums of self-adjoint free operators.  The main tool to
 study spectra of selfadjoint operators is the spectral measure and its
 Cauchy-Stieltjes transform. Recall that to every selfadjoint operator $T$ in a
 $C^*$-probability space with faithful state $\phi$ one can associate a
 probability measure $\mu$ which is characterized by the property that
 $\phi(T^n)=\int t^n d\mu(t)$. By faithfulness of the state $\phi$ one
 has $\supp\mu=\sigma(T)$. Thus knowing this measure one knows a
 fortiori the spectrum of $T$.  In order to compute this measure one
 has to solve a moment problem.  This is usually done via the Cauchy
 transform \eqref{eq:CauchyTransform}, which can be written in this
 case
 $$
 G(\zeta) = \int\frac{d\mu(t)}{\zeta-t}
 .
 $$
 This is an analytic function on  $\IC\setminus\sigma(T)$ mapping the
 upper half plane to the lower half plane.
 The absolute continuous part of the measure can be found by the
 Stieltjes inversion formula (cf.\
 \cite{Bercovici.Voiculescu:1998:regularity}) 
 $$
 \frac{d\mu}{dt}(t) = -\frac{1}{\pi} \lim_{\eps\searrow0} \Im G(t+i\eps)
 $$
 Atoms can be detected by studying the poles of $G(\zeta)$, namely
 $$
 \mu(\{t\}) = \lim_{\substack{\zeta\to t\\ \text{non-tangentially}}}
 (\zeta-t)\,G(\zeta)
 .
 $$

 When studying free sums of self-adjoint operators, by
 \eqref{eq:RTransformSum} one usually is given the inverse
 $K(z)=\frac{1}{z}(1 + \sum R_i(z))$ of the Cauchy transform $G(\zeta)$ and
 the latter is not accessible directly, but it is still possible to
 compute the spectrum. We are grateful to U.~Haagerup for showing the
 following lemma to us.

 %LLLLLLLLLLLLLLLLLLLLLLLLLLLLLLLLLLLLLLLLLLLLLLLLLLLLLLLLLLLLLLLLL
 \begin{Lemma}
   \label{lem:normal.spectrum.L2}
   Let $\alg{A}$ be a $C^*$-probability space with faithful state
   $\phi$, and let $L^2(\phi)$ be the associated $L^2$-space (the
   closure of $\alg{A}$ under the Hilbert norm
   $\norm{X}_2=\phi(X^*X)^{1/2}$). 
   Let $T$ be a normal operator, then
   $\sigma(T)=\overline{\{\lambda:\norm{(\lambda-T)^{-1}}_2=\infty\}}$.
 \end{Lemma}
 %lllllllllllllllllllllllllllllllllllllllllllllllllllllllllllllllll
 \begin{proof}
   Given $\lambda_0\in\sigma(T)$ and $\eps>0$ we have to find
   $\lambda\in\sigma(T)$ such that $\abs{\lambda-\lambda_0}<\eps$ and
   $\norm{(\lambda-T)^{-1}}_2=\infty$. To this end consider the square
   $Q_0$ of diameter $2\eps$ centered at $\lambda_0$. The spectral
   measure of this square is positive, say $\mu_T(Q_0)=\delta>0$. If
   follows that one of the four quarters of $Q_0$ has mass at least
   $\delta/4$. Let $Q_1$ be such a quarter. Doing the same argument
   again, we find a subsquare of $Q_1$ of mass at least
   $\delta/16$. Repeating this construction gives rise to a sequence
   $Q_n\subseteq Q_0$ with $\diam(Q_n)=\eps/2^{n-1}$ and mass
   $\mu_T(Q_n)\ge \delta/4^n$. Let $\lambda$ be the limit point of the
   sequence $Q_n$, then every point $\eta\in Q_n$ is at distance
   $\abs{\eta-\lambda}\leq \eps/2^n$ from $\lambda$ and thus 
   \begin{align*}
     \int\frac{d\mu_T(t)}{\abs{\lambda-t}^2}
     &\ge \frac{\mu_T(Q_0)}{\eps^2}
          + 
          \sum_{n=1}^\infty \mu_T(Q_n)
          \left(
            \frac{1}{(2^{-n}\eps)^2}
            -
            \frac{1}{(2^{1-n}\eps)^2}
          \right)\\
     &\ge \frac{\delta}{\eps^2}
          +
          \sum_{n=1}^\infty
           \frac{\delta}{4^n}
           \frac{4^n-4^{n-1}}{\eps^2}\\
     &=\frac{\delta}{\eps^2}
       \left(
         1+\sum_{n=1}^\infty \frac{3}{4}
       \right)\\
     &= \infty
   \end{align*}
 \end{proof}

 It is therefore enough to know the $L^2$-norm of the inverse, which in
 the case of a self-adjoint operator is easy to compute. For
 $\Im\zeta\ne0$ it is
 \begin{align*}
   \norm{(\zeta-T)^{-1}}_2^2
   &= \int\frac{d\mu_T(t)}{\abs{\zeta-t}^2} \\
   &= \int
       \frac{1}{\bar\zeta-\zeta}
       \left(
         \frac{1}{\zeta-t}
         -
         \frac{1}{\bar\zeta-t}
       \right)
       d\mu_T(t) \\
   &= - \frac{G(\zeta)-G(\bar\zeta)}{\zeta-\bar\zeta}
 \end{align*}
 while for $\Im\zeta=0$ we take the limit to get
 $$
 \norm{(\zeta-T)^{-1}}_2^2 = -G'(\zeta)
 .
 $$
 In terms of the inverse function $K(z)$ \eqref{eq:Ktransform} this
 becomes 
 %\QQQ[or Im $K(z)$?]
 $$
 \norm{(K(z)-T)^{-1}}_2^2 =
 \begin{cases}
   -\frac{z-\bar z}{K(z)-\overline{K(z)}} & \Im z\ne 0\\
   -\frac{1}{K'(z)}                       & \Im z = 0
 \end{cases}
 $$
 For concrete examples of the use of this see e.g.\ 
 \cite{Aomoto.Kato:1988:green,Figa-Talamanca.Steger:1994:harmonic}.

 %SSSSSSSSSSSSSSSSSSSSSSSSSSSSSSSSSSSSSSSSSSSSSSSSSSSSSSSSSSSSSSSSSSSSSS
 %SSSSSSSSSSSSSSSSSSSSSSSSSSSSSSSSSSSSSSSSSSSSSSSSSSSSSSSSSSSSSSSSSSSSSS
 \section{The non-selfadjoint case}
 \label{sec:nons.a.}
 %SSSSSSSSSSSSSSSSSSSSSSSSSSSSSSSSSSSSSSSSSSSSSSSSSSSSSSSSSSSSSSSSSSSSSS
 %SSSSSSSSSSSSSSSSSSSSSSSSSSSSSSSSSSSSSSSSSSSSSSSSSSSSSSSSSSSSSSSSSSSSSS

 Spectral theory of non-selfadjoint operators is complicated by the
 absence of the spectral measure\footnote{In some cases \emph{Brown's
     spectral measure} \cite{Brown:1986:lidskis} is computable and can
   serve a replacement -- see
   \cite{Haagerup.Larsen:1998:browns,Biane.Lehner:1999:computation}}
 and one has to analyze the resolvent directly.  We will see that the
 $L^2$-norm of the resolvent of free sums is easily computable and
 Haagerup inequality will provide a replacement of
 lemma~\ref{lem:normal.spectrum.L2}.

 %SSSSSSSSSSSSSSSSSSSSSSSSSSSSSSSSSSSSSSSSSSSSSSSSSSSSSSSSSSSSSSSSSSSSSS
 \subsection{The free resolvent}
 %SSSSSSSSSSSSSSSSSSSSSSSSSSSSSSSSSSSSSSSSSSSSSSSSSSSSSSSSSSSSSSSSSSSSSS

 The resolvent for the sum of convolution operators on free products of
 discrete groups has been computed by several authors, among the first
 are \cite{Woess:1986:nearest}, \cite{Steger:1985:harmonic},
 \cite{Soardi:1986:resolvent} and \cite{Cartwright.Soardi:1986:random}.
 On the free group itself it was also computed in
 \cite{Aomoto:1984:spectral} and \cite{Gerl.Woess:1986:local}.
 The resolvent of a sum of two free operators has been computed in
 \cite{Maassen:1992:addition} and very elegantly in
 \cite{Haagerup:1997:voiculescus}, see
 sec.~\ref{sec:TwoOperatorResolvent} below.

 Another proof for the free resolvent formula and the $R$-transform can
 be based on the following lemma (cf.\ 
 \cite{Figa-Talamanca.Steger:1994:harmonic},
 \cite{Lehner:1999:computing}).

 %LLLLLLLLLLLLLLLLLLLLLLLLLLLLLLLLLLLLLLLLLLLLLLLLLLLLLLLLLLLLLLLLLLLLLLLLL
 \begin{Lemma}
   \label{lem:FreeSummationFormula}
   Let $S_1,...,S_N\in B(H)$ be arbitrary operators and assume that the sum of
   alternating products
   $$
   \sum_{n=1}^\infty
    \sum_{i_1\ne i_2\ne \dots \ne i_n}
     S_{i_1}S_{i_2}\cdots S_{i_n}
   $$
   (the sum over all products where neighbouring factors are
   different) converges, then it equals
   $$
   \left(
     I - \sum_{i=1}^N S_i(I+S_i)^{-1}
   \right)^{-1}
   $$
 \end{Lemma}
 %lllllllllllllllllllllllllllllllllllllllllllllllllllllllllllllllllllllllll

 Convergence criteria will be presented below.
 The following is a simple reformulation of this summation formula.

 %PPPPPPPPPPPPPPPPPPPPPPPPPPPPPPPPPPPPPPPPPPPPPPPPPPPPPPPPPPPPPPPPPPPPPPPPP
 \begin{Proposition}
   \label{prop:GeneralResolventFormula}
   Let $X_1,\dots,X_N$ be operators on Hilbert space.
   Assume that $I-X_i$ is invertible for every
   $i\in\{1,\dots,N\}$ with inverse $(I-X_i)^{-1}=I+S_i$.
   Then
   $$
   \left(
     I - \sum X_i
   \right)^{-1}
   =I+\sum_{n=1}^\infty
       \sum_{i_1\ne i_2\ne\dots\ne i_n}
        S_{i_1}S_{i_2}\cdots S_{i_n},
   $$
   provided the sum on the right hand side converges.
 \end{Proposition}
 %ppppppppppppppppppppppppppppppppppppppppppppppppppppppppppppppppppppppppp

 Sums over alternating indices are well suited for free probability, as
 the very definition of freeness \eqref{eq:FreenessDefinition} suggests
 and we have the following theorem.

 %CCCCCCCCCCCCCCCCCCCCCCCCCCCCCCCCCCCCCCCCCCCCCCCCCCCCCCCCCCCCCCCCCCCCCCCCC
 \begin{Theorem}
   \label{thm:FreeResolvent}
   Let $(\alg{A},\phi)$ be a non-commutative probability space an let
   $T_i\in \alg{A}$ be freely independent random variables
   with Cauchy transforms
   $G_i(\zeta)=\phi( (\zeta-T_i)^{-1})$.
   The Cauchy transform is invertible at infinity with an inverse of
   the form $K_i(z)=\frac{1}{z}(1+R_i(z))$, i.e.\ 
   $\phi( (K_i(z)-T_i)^{-1} ) = z$ and thus
   $$
   (K_i(z)-T_i)^{-1} = z(1 + S_i(z))
   $$
   with $\phi(S_i(z))=0$.
   Then the resolvent of $T=\sum T_i$ at $K(z)=\frac{1}{z}(1+\sum R_i(z))$
   can be written formally as an infinite sum
   \begin{equation}
     \label{eq:FreeResolvent}
     (K(z) - \sum T_i)^{-1}
     = z \left(
           I
           +
           \sum_{n=1}^\infty \sum_{i_1\ne i_2\ne \dots\ne i_n}
            S_{i_1}S_{i_2}\cdots S_{i_n}
         \right)
   \end{equation}
   In particular, the expectation of the resolvent is
   $\phi
   \left(
     (K(z)-\sum T_i)^{-1}
   \right) =z
   $
   and \eqref{eq:RTransformSum} holds.

   If in addition
   \begin{equation}
     \label{eq:FreeResolventL2normBounded}
     \sum\frac{\norm{S_i}_2^2}{1+\norm{S_i}_2^2} < 1
   \end{equation}
   then the sum \eqref{eq:FreeResolvent} converges in $L^2(\phi)$ and its
   norm is
   \begin{equation}
     \label{eq:FreeResolventL2norm}
     \begin{aligned}
     \norm{
       (K(z)-\sum T_i)^{-1}
       }_2^2
     &=\abs{z}^2
       \left(
         1+ \sum_{n=1}^\infty \sum_{i_1\ne i_2\ne \dots\ne i_n}
         \norm{S_{i_1}}_2^2 \norm{S_{i_2}}_2^2 \cdots \norm{S_{i_n}}_2^2
       \right) \\
     &=\abs{z}^2
       \left(
         1-\sum\frac{\norm{S_i}_2^2}{1+\norm{S_i}_2^2}
       \right)^{-1}\\
     &=\abs{z}^2
       \left(
         1-\sum(1- \frac{\abs{z}^2}{\norm{(K_i(z)-T_i)^{-1}}_2^2})
       \right)^{-1}
     .
     \end{aligned}
   \end{equation}
   Moreover, boundedness of $\norm{(K(z)-\sum T_i)^{-1}}_2$ implies
   boundedness of $\norm{(K(z)-\sum T_i)^{-1}}$ (operator norm).
 \end{Theorem}
 %ccccccccccccccccccccccccccccccccccccccccccccccccccccccccccccccccccccccccc

 %Note that formula \eqref{eq:RTransformSum} for the $R$-transform is an
 %easy consequence of the preceding theorem.

 The reader easily verifies the following reformulation
 in terms of moment generating functions.

 \begin{Corollary}
   Let $a_i$, $i=1,\dots,n$ be freely independent random variables in the
   noncommutative probability space $(\alg{A},\phi)$.
   Let $f_i(s)=\frac{1}{s} G_i(\frac{1}{s}) = \phi ((1-sa_i)^{-1})$ the
   moment generating function and denote
   $$
   \overset\circ{a}_i(s) = (1-sa_i)^{-1}-f_i(s)
   $$
   the centered part of the resolvent.
   Given $\lambda\in\IC$, assume that
   there exist numbers $z\ne0$ and $s_i$ such that
   $s_if_i(s_i)=z$ for each $i=1,\dots,n$ and such that
   $\lambda=\frac{1}{z}+\sum(\frac{1}{s_i}-\frac{1}{z})$
   \begin{enumerate}
    \item If
     \begin{equation}
       \label{eq:FreeResolventL2normBoundedGF}
       \sum \frac{\snorm{\overset\circ{a}_i}_2^2}{1+\snorm{\overset\circ{a}_i}_2^2}<1
     \end{equation}
     then $\lambda\not\in\sigma(\sum a_i)$.
    \item If
     $$
     \sum \frac{\snorm{\overset\circ{a}_i}_2^2}{1+\snorm{\overset\circ{a}_i}_2^2}=1
     $$
     and there are $\lambda'$ arbitrary close to $\lambda$ for which
     \eqref{eq:FreeResolventL2normBoundedGF} holds,
     then $\lambda\in\sigma(\sum a_i)$.
   \end{enumerate}
 \end{Corollary}

 Theorem~\ref{thm:FreeResolvent} provides an easy to check sufficient
 criterion for the boundedness in $L^2(\phi)$.  By virtue of the
 following variant of the Haagerup inequality we will see that it is
 also sufficient for boundedness in operator norm.

 %PPPPPPPPPPPPPPPPPPPPPPPPPPPPPPPPPPPPPPPPPPPPPPPPPPPPPPPPPPPPPPPPPPPPPPPPP
 \begin{Proposition}[{\cite[Lemma~3.4]{Dykema.Haagerup.Rordam:1997:stable}}]
   \label{prop:HaagerupInequality}
   Let $(\alg{A},\tau)$ be a tracial $C^*$-probability space and let
   $\alg{A}_i$ be free subalgebras with %standard\footnote{define this}
   orthonormal bases $X_i=\{I,x_{i,1},x_{i,2},...\}$. Let
   $\overset\circ{X}_i= X_i\setminus \{I\}$ be the centered part and
   $$
   X = \{x_{1,k_1}x_{2,k_2}\dots x_{n,k_n} : x_{j,k_j}\in
       \overset\circ{X}_{i_j}, i_j\ne i_{j+1}\}\cup \{I\}
   $$
   the free product of these bases.  This is an orthonormal basis of
   the free product $\ast \alg{A}_i$ and can be decomposed $X=\bigcup_n
   E_n$ into the subsets $E_n$ of words of length $n$.  For a finitely
   supported operator $a\in\linspan X$ denote by $F_i(a)$ its
   $i$-support, i.e., the set of all $x\in \overset\circ{X}_i$
   appearing in the words of the expansion of $a$.  Then for any $a\in
   \linspan E_n$ we have
   $$
   \norm{a} \le (2n+1)
            \max_i \left(
                     \sum_{x\in F_i(a)} \norm{x}^2
                   \right)^{1/2}
            \norm{a}_2
   $$
 \end{Proposition}
 %ppppppppppppppppppppppppppppppppppppppppppppppppppppppppppppppppppppppppp

 \begin{proof}[Proof of Theorem~\ref{thm:FreeResolvent}]
   Formula \eqref{eq:FreeResolvent} follows from
   Proposition~\ref{prop:GeneralResolventFormula}.
   For the $L^2$-norm, note that the freeness condition implies that 
   for different sets of indices
   $i_1,i_2,\dots,i_m$ and $j_1,j_2,\dots,j_n$ the summands
   $S_{i_1}S_{i_2}\dots S_{i_m}$ and  $S_{j_1}S_{j_2}\dots S_{j_n}$ are
   orthogonal and that
   $\norm{S_{i_1}S_{i_2}\dots S_{i_m}}_2^2
   =\norm{S_{i_1}}_2^2\norm{S_{i_2}}_2^2\dots\norm{S_{i_m}}_2^2$.
   Then formula \eqref{eq:FreeResolventL2norm} follows from
   Lemma~\ref{lem:FreeSummationFormula}. To show boundedness of the
   resolvent in operator norm, we resort to analytic functions, as in
   \cite[Chapter~2, Lemma~1.6]{Figa-Talamanca.Steger:1994:harmonic}.
   Assume that the $L^2$-norm is bounded, i.e.,
   $\sum \norm{S_{i_1}}_2^2\dots \norm{S_{i_n}}_2^2<\infty$.
   Then we can define an analytic function on the open disk $\ID$
   $$
   F(\xi) = 1 +
            \sum_{n=1}^\infty
            \sum_{i_1\ne i_2\ne\dots\ne i_n}
             \xi^n
             \norm{S_{i_1}}_2^2\dots \norm{S_{i_n}}_2^2
   $$
   which by the summation formula from Lemma~\ref{lem:FreeSummationFormula}
   equals
   $$
   F(\xi)=
   \left(1-
     \sum \frac{\xi \norm{S_i}_2^2}{1+\xi \norm{S_i}_2^2}
   \right)^{-1}
   .
   $$
   This is a rational function and by assumption it has no pole on
   the circle $\{\xi:\abs{\xi}=1\}$. As a rational function it has 
   finitely many singularities and is therefore analytic on some disk of radius
   $1+\eps$ with $\eps>0$.  It follows that there exists some constant
   $C$ independent of $n$ such that the Taylor coefficients satisfy
   $$
   \sum_{i_1\ne i_2\ne\dots\ne i_n}
    \norm{S_{i_1}}_2^2\dots \norm{S_{i_n}}_2^2
    \le C (1+\eps)^{-n}
   $$
   Now with the help of the Haagerup inequality
   (Proposition~\ref{prop:HaagerupInequality}) we obtain the estimate
   \begin{align*}
     \Bignorm{
       I
       +
       \sum_{n=1}^\infty \sum_{i_1\ne i_2\ne \dots\ne i_n}
        S_{i_1}S_{i_2}\cdots S_{i_n}
     } \\
     \phantom{=}&\le 1 + 
          \sum_{n=1}^\infty
           (2n+1) \max\norm{S_i} 
           \left(
             \sum_{i_1\ne i_2\ne\dots\ne i_n}
             \norm{S_{i_1}}_2^2\dots \norm{S_{i_n}}_2^2
           \right)^{1/2} \\
     &\le 1 + 
          \sum_{n=1}^\infty
           (2n+1) \max\norm{S_i} 
           \sqrt{C} (1+\eps)^{-n/2} \\
     &<\infty
   \end{align*}
 \end{proof}

 The criterion for boundedness of the $L^2$-norm
 \begin{equation}
   \label{eq:sumSi2}
   \sum\frac{\norm{S_i}_2^2}{1+\norm{S_i}_2^2}
   < 1
 \end{equation}
 can be reformulated in terms of symmetric functions.
 Recall that the \emph{elementary symmetric functions} in $n$ variables
 $x_1,x_2,\dots,x_n$ are defined as 
 $$
 E_k = E_k(x_1,x_2,\dots,x_n)
     = \sum_{i_1<\dots< i_k}
        x_{i_1}x_{i_2}\dots x_{i_k}
 ;
 $$
 their generating function is
 $$
 \prod_{i=1}^n (1+tx_i) = 1 + \sum_{k=1}^n t^k E_k
 $$
 The left hand side of \eqref{eq:sumSi2} is a rational symmetric function of
 $x_i=\norm{S_i}_2^2$ and can be expressed in terms of the elementary
 symmetric functions as follows.
 \begin{align*}
   \sum\frac{x_i}{1+t x_i}
   &= \frac{\sum_i x_i\prod_{j\ne i}(1+tx_j)}{\prod(1+tx_i)} \\
   &= \frac{\frac{d}{dt}\prod(1+tx_j)}{\prod(1+tx_i)} \\
   &= \frac{\sum_{k=1}^n k\, E_k t^{k-1}}{1+\sum_{k=1}^n  E_k t^k}
 \end{align*}
 at $t=1$ and $x_i=\norm{S_i}_2^2$ we get therefore the condition
 $$
 \frac{\sum_{k=1}^n k\, E_k(\norm{S_1}_2^2,\dots, \norm{S_n}_2^2)}%
      {1+\sum_{k=1}^n E_k(\norm{S_1}_2^2,\dots,\norm{S_n}_2^2)}<1
 $$
 which is equivalent to
 $$
 \sum_{k=1}^n (k-1)\, E_k(\norm{S_1}_2^2,\dots,\norm{S_n}_2^2) <1
 $$
 For $n=2$ this reduces to the simple condition
 $\norm{S_1}_2^2\, \norm{S_2}_2^2<1$. We will investigate this case in
 section~\ref{sec:TwoOperatorResolvent} below.

 %SSSSSSSSSSSSSSSSSSSSSSSSSSSSSSSSSSSSSSSSSSSSSSSSSSSSSSSSSSSSSSSSSSSSSS
 %SSSSSSSSSSSSSSSSSSSSSSSSSSSSSSSSSSSSSSSSSSSSSSSSSSSSSSSSSSSSSSSSSSSSSS
 \section{Case of two operators}
 \label{sec:TwoOperatorResolvent}
 %SSSSSSSSSSSSSSSSSSSSSSSSSSSSSSSSSSSSSSSSSSSSSSSSSSSSSSSSSSSSSSSSSSSSSS
 %SSSSSSSSSSSSSSSSSSSSSSSSSSSSSSSSSSSSSSSSSSSSSSSSSSSSSSSSSSSSSSSSSSSSSS
 In the particular case of the sum of two operators the analysis can be
 somewhat simplified and there are shorter proofs.  We can use the
 following proposition instead of the Haagerup inequality from
 Proposition~\ref{prop:HaagerupInequality}. It is much easier to prove.

 %PPPPPPPPPPPPPPPPPPPPPPPPPPPPPPPPPPPPPPPPPPPPPPPPPPPPPPPPPPPPPPPPPPPPPPPPP
 \begin{Proposition}[{\cite[Prop.~4.1]{Haagerup.Larsen:1998:browns}}]
   \label{prop:HaagerupLarsen:SpectralRadius}
   Let $(\alg{M},\tau)$ be a non-commutative von Neumann probability space with faithful
   trace state $\tau$ and let $a$, $b\in\alg{M}$ be arbitrary centered
   ($\tau(a)=\tau(b)=0$) $*$-free random variables. Then the spectral
   radius of their product is
   $$
   \rho(ab) = \norm{ab}_2 = \norm{a}_2 \norm{b}_2
   $$
 \end{Proposition}
 %ppppppppppppppppppppppppppppppppppppppppppppppppppppppppppppppppppppppppp

 %CCCCCCCCCCCCCCCCCCCCCCCCCCCCCCCCCCCCCCCCCCCCCCCCCCCCCCCCCCCCCCCCCCCCCCCCC
 \begin{Corollary}
   \label{cor:RadiusabContainsT}
   For $a$ and $b$ as in
   Proposition~\ref{prop:HaagerupLarsen:SpectralRadius}, the circle of
   radius $\rho(ab)=\norm{ab}_2$ is part of the spectrum $\sigma(ab)$.
 \end{Corollary}
 %ccccccccccccccccccccccccccccccccccccccccccccccccccccccccccccccccccccccccc

 \begin{proof}
   We can rescale $a$ and $b$ and assume that
   $\norm{a}_2=\norm{b}_2=1$. We have to show that $1-tab$ is not
   invertible whenever $\abs{t}=1$.  whenever $\abs{t}<1$, then $1-tab$
   is invertible and $\frac{1}{t}$ is in the complement of the
   spectrum.  By orthogonality, the $L^2$-norm of the inverse is
   $$
   \norm{(1-tab)^{-1}}_2^2
   = \sum_{n=0}^\infty
      \abs{t}^{2n}\,
      \norm{a}_2^{2n}\,
      \norm{b}_2^{2n} 
   = \frac{1}{1-\abs{t}^2\, \norm{a}_2^2\, \norm{b}_2^2}
   $$
   and this grows unboundedly, as $\abs{t}$ tends to one.
   Consequently the operator norm of the inverse becomes unbounded as
   $\abs{t}\to 1$. Since the resolvent is continuous (even analytic) on
   the complement of the spectrum, any number $t$ of modulus $\abs{t}=1$
   is in the spectrum of $ab$.
 \end{proof}

 Part \eqref{prop:2ResolventFormula:i} of the following proposition is taken from
 {\cite{Haagerup:1997:voiculescus}}.

 %PPPPPPPPPPPPPPPPPPPPPPPPPPPPPPPPPPPPPPPPPPPPPPPPPPPPPPPPPPPPPPPPPPPPPPPPP
 \begin{Proposition}
   \label{prop:2ResolventFormula}
   Let $a$, $b\in\alg{A}$ and $\abs{s}<\frac{1}{\rho(a)}$,
   $\abs{t}<\frac{1}{\rho(b)}$,
   $f(s)=\phi((1-sa)^{-1})$,
   $g(t)=\phi((1-tb)^{-1})$.
   Put $\overset{\circ}{a}(s)=(1-sa)^{-1}-f(s)$,
   $\overset{\circ}{b}(t)=(1-tb)^{-1}-g(t)$ and
   assume that 
   \begin{equation}
     \label{eq:TwoOp:sfs=tgt}
     sf(s)=tg(t) \ne 0
     .
   \end{equation}
   Set
   \begin{equation}
     \label{eq:TwoOp:lambda=1s+1t-1sfs}
     \lambda =
     \frac{f(s)+g(t)-1}{sf(s)}=\frac{1}{s}+\frac{1}{t}-\frac{1}{sf(s)},
   \end{equation}
   then
   \begin{enumerate}[(i)]
    \item
     \label{prop:2ResolventFormula:i}
     $$
     \lambda-a-b
     = \frac{g(t)}{s}
       (1-sa)
       \left(
         1 - \frac{\overset{\circ}{a}(s)\,
                   \overset{\circ}{b}(t)}%
                  {f(s)\, g(t)}
       \right)
       (1-tb)
     .
     $$
    \item
     \label{prop:2ResolventFormula:ii}
     If 
     $
     \snorm{\overset{\circ}{a}(s)}_2^2\,
     \snorm{\overset{\circ}{b}(t)}_2^2 < \abs{f(s)}^2\,
     \abs{g(t)}^2$,
     then
     $\lambda\in \IC\setminus\sigma(a+b)$.
    \item
     \label{prop:2ResolventFormula:iii}
     If
     $\snorm{\overset{\circ}{a}(s)}_2^2\,
     \snorm{\overset{\circ}{b}(t)}_2^2=\abs{f(s)}^2\,\abs{g(t)}^2$,
     then
     $\lambda\in\sigma(a+b)$.
   \end{enumerate}
 \end{Proposition}
 %ppppppppppppppppppppppppppppppppppppppppppppppppppppppppppppppppppppppppp

 \begin{proof}
   Part~\eqref{prop:2ResolventFormula:i} can be derived from
   Theorem~\ref{thm:FreeResolvent} and can also be verified directly by
   expanding the right hand side.
   Parts~\eqref{prop:2ResolventFormula:ii} and
   \eqref{prop:2ResolventFormula:iii} follow from
   Proposition~\ref{prop:HaagerupLarsen:SpectralRadius} and
   Corollary~\ref{cor:RadiusabContainsT} respectively and the fact that
   $1-sa$ and $1-tb$ are invertible.
 \end{proof}

 %SSSSSSSSSSSSSSSSSSSSSSSSSSSSSSSSSSSSSSSSSSSSSSSSSSSSSSSSSSSSSSSSSSSSSS
 %SSSSSSSSSSSSSSSSSSSSSSSSSSSSSSSSSSSSSSSSSSSSSSSSSSSSSSSSSSSSSSSSSSSSSS
 \section{Examples}
 \label{sec:Examples}
 %SSSSSSSSSSSSSSSSSSSSSSSSSSSSSSSSSSSSSSSSSSSSSSSSSSSSSSSSSSSSSSSSSSSSSS
 %SSSSSSSSSSSSSSSSSSSSSSSSSSSSSSSSSSSSSSSSSSSSSSSSSSSSSSSSSSSSSSSSSSSSSS

 The most interesting examples are perhaps those coming from
 convolution operators on free products of discrete groups.  Let us
 consider here non-selfadjoint convolution operators supported on the
 generators of free products of cyclic groups.  For a survey on spectra
 of such operators on more general finitely generated discrete groups
 see
 \cite{Harpe.Robertson.Valette:1993:spectrum,Harpe.Robertson.Valette:1993:spectrum*1}.
 Working with finite groups has the advantage that the involved
 equations are algebraic and the powerful machinery of algebraic
 geometry is available to obtain quite explicit results.  Let us
 therefore recall some facts from algebraic geometry.  For more
 background on algebraic equations we refer to the excellent survey
 \cite{Cox.Little.OShea:1998:using} and text book
 \cite{Cox.Little.OShea:1997:ideals}.

 We used \verb|Singular|
 \cite{Greuel.Pfister.Schoenemann:1998:singular-version-1-2-user-manual-}
 for the algebraic computations, \verb|GNU octave| \cite{Eaton:1998:gnu} for
 the eigenvalue computations and \verb|Mathematica|
 \cite{Wolfram:1996:mathematica} for visualizations.

 \subsection{Eigenvalue approach to polynomial equations}
 \label{sub:ex:eigenvalueapproach}

 The most efficient numerical approach seems to be the 
 \emph{matrix eigenvalue method} following
 \cite[section~2.4]{Cox.Little.OShea:1998:using}.  Recall that systems
 of polynomial equations are in one-to-one correspondence with
 polynomial ideals. Let $f_1,\dots,f_m\in\IC[x_1,\dots,x_n]$ be
 polynomials and denote by $V=V(f_1,\dots,f_m)$ the set of solutions
 (also called \emph{algebraic set} or \emph{algebraic variety}) of the
 system of polynomial equations
 \begin{equation}
   \label{eq:polynomialsystem}
   f_1=0, \dots, f_m=0.
 \end{equation}
 On the other hand, on can associate to this system the ideal
 $$
 I=\langle f_1,\dots,f_m \rangle
  =\left\{
     \sum g_i f_i : g_i\in \IC[x_1,\dots,x_n]
   \right\}
 $$
 it generates. Then it is easy to see that the variety $V$ does not
 depend on the generators, but only on the ideal $I$. That is,
 $V=V(I) = \{x\in \IC: f(x)=0\ \forall f\in I\}$
 (note however, that different ideals can lead to
 the same variety $V$).

 There are close relations between the Properties of the ideal and its
 variety.  In particular, one can show that the system
 \eqref{eq:polynomialsystem} has finitely many solutions ($V$ is
 \emph{zero-dimensional}) if and only if the quotient algebra
 $A=\IC[x_1,x_2,\dots,x_n]/I$ is finite dimensional.  Then for any
 polynomial $h\in\IC[x_1,\dots,x_n]$, the eigenvalues of multiplication
 map
 \begin{align*}
   M_h:A&\rightarrow A\\
   [g]&\mapsto [hg]
 \end{align*}
 coincide with the values of the polynomial $h(\xi)$, evaluated at
 $\xi\in V$.  In particular, choosing $h=x_i$ one can compute the
 coordinates of the solutions; in fact the eigenvalues and
 eigenvectors of a single matrix $M_h$ are sufficient to determine all
 the solutions, see \cite[section~2.4]{Cox.Little.OShea:1998:using}.
 Such an approach avoids the propagation of rounding errors which are
 common to elimination methods.

 In order to do concrete computations in the quotiont algebra $A$ one
 needs a normal forms for its elements. This can be accomplished using
 Gr\"obner bases.
 Given any (linear) order on the monomials, one can use a generalized
 Euclid's algorithm to write any $f\in \IC[x_1,\dots,x_n]$ as a sum
 $$
 f = q_1 f_1 + \dots + q_m f_m + r
 $$
 The remainder $r$ of this division is a natural candidate for a
 normal form of $f$ modulo the ideal $I$, however it may not be unique
 and there may exist $f\in I$ with $r\ne0$. Fortunately there always
 exist generating sets which are nice in this respect.  A
 \emph{Gr\"obner basis} of the ideal $I$ is a generating set
 $G=\{g_1,\dots,g_p\}$ s.t.\ the remainder $r$ of Euclidean division is
 unique, i.e., $f\in I$ if and only if $r=0$.  Gr\"obner bases always
 exist and can be constructed with Buchberger's algorithm. With such
 bases the remainder $r$ can serve as a normal form for the elements of
 $A$.  Another characterizing property of Gr\"obner bases is the
 following.  For a polynomial $f\in \IC[x_1,\dots,x_n]$ denote $L(f)$
 the leading monomial of $f$ (with respect to the chosen monomial
 order) and let $L(I)$ be the ideal generated by the leading monomials
 of all elements of $I$. Then $G$ is a Gr\"obner basis if and only if
 $\langle L(G)\rangle = L(I)$.  In the case of a zero dimensional
 ideal, a natural basis of $A=\IC[x_1,\dots,x_n]/I$ are the equivalence
 classes
 $$
 B = \{[x^\alpha] : x^\alpha \not\in \langle L(G)\rangle\}
 $$
 i.e.\ the (finitely many) equivalence classes of monomials which
 are reduced with respect to the Gr\"obner basis $G$.
 %
 %to compute a reduced form 
 %compute the matrix of the multiplication map $M_{x_i}$,
 %one has to work with a basis of $A$, which can be chosen to be the
 %equivalence classes of monomials
 %$B=\{[x^k]=[x_1^{k_1}\dots x_n^{k_n}] : k\in K\}$
 %where $K$ is a finite set of multiexponents, adapted to some
 %\emph{Gr\"obner basis} of the ideal $I$ with respect to some monomial
 %ordering. 
 %
 %Recall that a Gr\"obner basis of an ideal
 %is a generating set of polynomials with respect to which the
 %generalized division algorithm yields a unique remainder, which can be
 %used as a normal form for the elements of the quotient algebra $A$.
 %Equivalently, a Gr\"obner basis $G$ is a generating set for which 
 %the ideal generated by the leading terms $\langle LT(G)\rangle=\langle
 %LT(I)\rangle$.

 \subsection{Elimination}
 \label{sub:ex:elimination}

 Gr\"obner bases are also needed for elimination.
 Roughly speaking, elimination of variables corresponds to projection
 of algebraic varieties.
 Elimination is done by computing the intersection $I_k=I\cap
 \IC[x_{k+1},\dots,x_n]$. For suitable monomial orders (e.g.\
 lexicographical order), one can use the fact that a Gr\"obner basis
 $G$ of $I$ has the property that $G_k=G\cap \IC[x_{k+1},\dots,x_n]$ is a
 Gr\"obner basis for $I_k$. 
 However, lexicographical Gr\"obner bases
 are expensive to compute and it is usually more efficient to work with
 other monomial orders, like \emph{degree reverse lexicographical
   order} and use other methods for elimination.

 \subsection{Free product groups}
 \label{sub:ex:freeproductgroups}

 For simplicity, let us consider the sum of two convolution operators
 $u_m+v_n$, where $u_m,v_n$ are the generators of
 $C_\lambda^*(\IZ_m*\IZ_n)$. 
 The moments of $u_m$ are
 $$
 \tau(u_m^k) = 
 \begin{cases}
   1 & k = 0 \mod m\\
   0 & \text{otherwise}
 \end{cases}
 $$
 and the inverse of $1-su_m$ is
 $$
 (1-su_m)^{-1} = \frac{1+su_m+\dots+s^{m-1}u_m^{m-1}}{1-s^m}
 .
 $$
 Thus for the moment generating function we obtain
 $$
 f(s) = \tau((1-su_m)^{-1}) = \frac{1}{1-s^m}
$$
From Proposition~\ref{prop:2ResolventFormula} we infer that if
the equations \eqref{eq:TwoOp:sfs=tgt} and
\eqref{eq:TwoOp:lambda=1s+1t-1sfs} satisfied, which in our particular example
read
\begin{gather*}
  \frac{s}{1-s^m} = \frac{t}{1-t^n} \\
  \lambda = \frac{1}{t} + s^{m-1}
\end{gather*}
and if $s$ and $t$ are such that the $L^2$-norm condition
$$
\frac{\abs{s}^2+\abs{s}^4+\dots+\abs{s}^{2m-2}}%
{\abs{1-s^m}^2}
\,
\frac{\abs{t}^2+\abs{t}^4+\dots+\abs{t}^{2n-2}}%
{\abs{1-t^n}^2}
<\frac{1}{\abs{1-s^m}^2}\, \frac{1}{\abs{1-t^n}^2}
$$
holds, then $\lambda$ is not in the spectrum of $u_m+u_n$.  In
order to find the outer border of the spectrum (that is, the curve
where the $L^2$-norm becomes infinite), we have to solve the equations
\begin{subequations}
  \label{eq:um+un:AlgEqu}
  \begin{gather}
    \label{eq:um+un:AlgEqu:a}
    s(1-t^n)-t(1-s^m) = 0\\
    \label{eq:um+un:AlgEqu:b}
    \lambda t - 1 - s^{m-1}t = 0\\
    \label{eq:um+un:AlgEqu:c}
    (\abs{s}^2+\abs{s}^4+\dots+\abs{s}^{2m-2})\,
    (\abs{t}^2+\abs{t}^4+\dots+\abs{t}^{2n-2})<1
  \end{gather}
\end{subequations}
and look where the inequality \eqref{eq:um+un:AlgEqu:c} becomes an equality.
There are two approaches to cope with the non-algebraic part \eqref{eq:um+un:AlgEqu:c}.
One is to solve the algebraic part
\eqref{eq:um+un:AlgEqu:a}, \eqref{eq:um+un:AlgEqu:b} first and verify
\eqref{eq:um+un:AlgEqu:c} numerically.
The alternative approach is to introduce new variables for the real
and imaginary parts of $\lambda=x+iy$, $s=s_1+is_2$ and $t=t_1+it_2$,
and to separate real and imaginary parts of the equations, which now
become purely algebraic, but with more unknowns. The real solutions of
the new system correspond to the complex solutions of the original
one.

The second approach makes sense if the larger system is sufficiently
small to allow complete algebraic elimination of the parameters $s_i$
and $t_i$; for a numerical solution the first approach is preferable.

\begin{figure}[t]
  \label{fig:u2+u3}
  \begin{center}
  \includegraphics[width=\textwidth, height=.4\textheight, keepaspectratio]{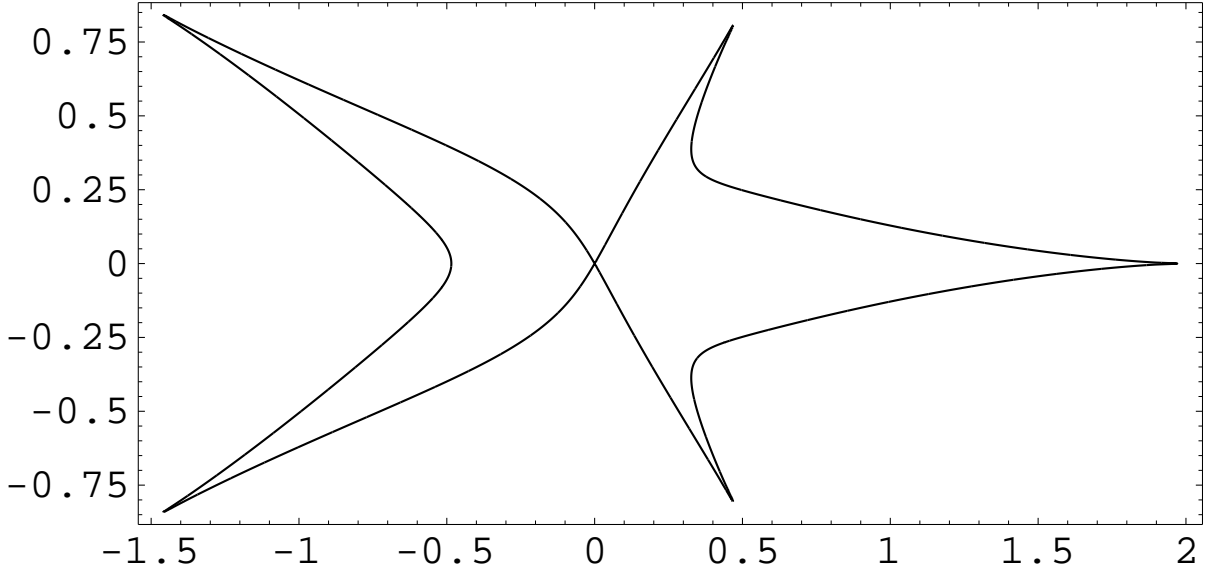}
    \caption{Spectrum of $u_2+u_3$%, $1050\times 900$ points
}
  \end{center}
\end{figure}

\begin{Example}[``$u_2+v_3$'']
Let us consider the simplest non-trivial example ($m=2$ and $n=3$) and
compute the border of the spectrum of $u_2+v_3$. If for a given
$\lambda$ the system
\begin{subequations}
  \label{eq:u2+u3:AlgEqu}
  \begin{gather}
    \label{eq:u2+u3:AlgEqu:a}
    s(1-t^3)-t(1-s^2) = 0\\
    \label{eq:u2+u3:AlgEqu:b}
    \lambda t - 1 - st = 0\\
    \label{eq:u2+u3:AlgEqu:c}
    \abs{s}^2 (\abs{t}^2+\abs{t}^4)<1
  \end{gather}
\end{subequations}
has a solution, then $\lambda$ is not in the spectrum of $u_2+v_3$.

In order to get the border of the spectrum, we replace the inequality
\eqref{eq:u2+u3:AlgEqu:c} by an equality and
separate real and imaginary part, which leads to the system

\begin{equation}
  \label{eq:u2+u3:SpectralBorderSystem}
\begin{aligned}
-s_{1}\,t_{1}^{3}+3\,s_{2}\,t_{1}^{2}\,t_{2}+3\,s_{1}\,t_{1}\,t_{2}^{2}-s_{2}\,t_{2}^{3}
+s_{1}^{2}\,t_{1}-s_{2}^{2}\,t_{1}-2\,s_{1}\,s_{2}\,t_{2}+s_{1}-t_{1} &=0 \\
-s_{2}\,t_{1}^{3}-3\,s_{1}\,t_{1}^{2}\,t_{2}+3\,s_{2}\,t_{1}\,t_{2}^{2}+s_{1}\,t_{2}^{3}
+2\,s_{1}\,s_{2}\,t_{1}+s_{1}^{2}\,t_{2}-s_{2}^{2}\,t_{2}+s_{2}-t_{2} &=0 \\
-x\,t_{1}+s_{1}\,t_{1}+y\,t_{2}-s_{2}\,t_{2}+1&= 0\\
-y\,t_{1}+s_{2}\,t_{1}-x\,t_{2}+s_{1}\,t_{2} &= 0\\
s_{1}^{2}\,t_{1}^{4}+s_{2}^{2}\,t_{1}^{4}+2\,s_{1}^{2}\,t_{1}^{2}\,t_{2}^{2}+2\,s_{2}^{2}\,t_{1}^{2}\,t_{2}^{2}+s_{1}^{2}\,t_{2}^{4}+s_{2}^{2}\,t_{2}^{4}+s_{1}^{2}\,t_{1}^{2}+s_{2}^{2}\,t_{1}^{2}+s_{1}^{2}\,t_{2}^{2}+s_{2}^{2}\,t_{2}^{2}-1 &=0
.
\end{aligned}
\end{equation}
Using degrevlex order there is a Gr\"obner basis of size 34 for the
ideal $I$ generated by the left hand sides of these equations and elimination
of $s_i$ and $t_i$ succeeds: the ideal $I\cap \IC[x,y]$ is generated
by a single polynomial and we get the implicit equation
\begin{gather*}
  x^{16}+8\,x^{14}\,y^2+28\,x^{12}\,y^4+56\,x^{10}\,y^6+70\,x^8\,y^8+56\,x^6\,y^{10} +28\,x^4\,y^{12}+8\,x^2\,y^{14}+y^{16}\\
 -6\,x^{15}+54\,x^{13}\,y^2+226\,x^{11}\,y^4
 +238\,x^9\,y^6-18\,x^7\,y^8-158\,x^5\,y^{10}-74\,x^3\,y^{12}-6\,x\,y^{14}\\
+5\,x^{14}
-65\,x^{12}\,y^2+305\,x^{10}\,y^4+435\,x^8\,y^6-225\,x^6\,y^8-179\,x^4\,y^{10}
+107\,x^2\,y^{12}+y^{14}\\
+32\,x^{13}-400\,x^{11}\,y^2-400\,x^9\,y^4+480\,x^7\,y^6
+384\,x^5\,y^8-80\,x^3\,y^{10}-16\,x\,y^{12}\\
-59\,x^{12}+74\,x^{10}\,y^2-1033\,x^8\,y^4
-548\,x^6\,y^6+547\,x^4\,y^8-70\,x^2\,y^{10}+y^{12}\\
-40\,x^{11}+776\,x^9\,y^2 -1008\,x^7\,y^4+112\,x^5\,y^6-104\,x^3\,y^8
+ 8\,x\,y^{10}\\
+136\,x^{10}+48\,x^8\,y^2 +736\,x^6\,y^4-176\,x^4\,y^6+24\,x^2\,y^8\\
-32\,x^9+224\,x^7\,y^2-224\,x^5\,y^4 +32\,x^3\,y^6\\
-48\,x^8-32\,x^6\,y^2+16\,x^4\,y^4   = 0
.
\end{gather*}
This implicit equation defines the curve shown in figure~\ref{fig:u2+u3}.
The figure does not show the isolated solutions $2$,
$\frac{1}{2}\pm\frac{\sqrt{3}}{2}i$,
$-\frac{3}{2}\pm\frac{\sqrt{3}}{2}i$, which however are not in the
spectrum, as they also come from solutions of
(\ref{eq:u2+u3:AlgEqu:a}, \ref{eq:u2+u3:AlgEqu:b})
which satisfy \eqref{eq:u2+u3:AlgEqu:c}.

The intersection with the $x$-axis are the solutions at $y=0$:
\begin{multline*}
  x^{16}-6\,x^{15}+5\,x^{14}+32\,x^{13}-59\,x^{12}-40\,x^{11}+136\,x^{10}-32\,x^9-48\,x^8 \\
= (x^4+2\,x^3-3\,x^2-8\,x-3)(x-2)^4 x^8=0
\end{multline*}
in particular, the spectral radius is the positive solution of
$$
x^4+2\,x^3-3\,x^2-8\,x-3=0
$$
which is
\begin{equation}
  \label{eq:u2+u3:Spectralradius}
  \begin{aligned}
    \rho(u_2+u_3)
% &=\frac{1}{6}
%   \Biggl(
%     -3 + {\sqrt{3\,\left( 9 + {{\left( 135 - 6\,{\sqrt{249}} \right) }^{{\frac{1}{3}}}} + 
%            {{\left( 135 + 6\,{\sqrt{249}} \right) }^{{\frac{1}{3}}}}
%            \right) }} \\
%   &\phantom=\left.
%      + 3\,{\sqrt{6 - {\frac{{{\left( 135 - 6\,{\sqrt{249}} \right) }^{{\frac{1}{3}}}}}
%             {3}} - {\frac{{{\left( 135 + 6\,{\sqrt{249}} \right) }^{{\frac{1}{3}}}}}{3}} + 
%           8\,{\sqrt{{\frac{3}
%                 {9 + {{\left( 135 - 6\,{\sqrt{249}} \right) }^{{\frac{1}{3}}}} + 
%                   {{\left( 135 + 6\,{\sqrt{249}} \right) }^{{\frac{1}{3}}}}
%                   }}}}}}
%   \;\right) \\
    &=\frac{1}{6}
    \Biggl\{
      -3
      + \sqrt{3\,
        \left( 9 + \left(135 - 6\,\sqrt{249} \right)^{\frac{1}{3}}
          + 
          \left( 135 + 6\,{\sqrt{249}} \right) ^{\frac{1}{3}}
        \right) } \\
    &\phantom=
    + 3 
      \left(
        6
        -
        \frac{\left( 135 - 6\,\sqrt{249} \right)^{\frac{1}{3}}}%
             {3}
        - \frac{\left( 135 + 6\,{\sqrt{249}}\right)^{\frac{1}{3}}}%
               {3}
      \right. \\
    &\phantom=\qquad
    + \left.
        8\,\sqrt{{\frac{3}%
                       {9 + {{\left( 135 - 6\,{\sqrt{249}} \right) }^{{\frac{1}{3}}}} + 
                {{\left( 135 + 6\,{\sqrt{249}} \right) }^{{\frac{1}{3}}}}
                }}}
      \right)^{\frac{1}{3}}
    \;\Biggr\} \\
 &\simeq 1.97148
  \end{aligned}
\end{equation}
With the exception of the double point at $0$ the apparent
singularities of the curve only appear as such and the curve is
actually smooth; for example, the curvature at the point of maximal
modulus calculated in \eqref{eq:u2+u3:Spectralradius} is approximately
$166053.0$, so the radius of curvature is about $1/166053.0\simeq
6.02219 \times 10^{-6}$.
\end{Example}

\begin{figure}[t]
  \label{fig:uuivv}
  \begin{center}
  \includegraphics[width=\textwidth, height=.4\textheight, keepaspectratio]{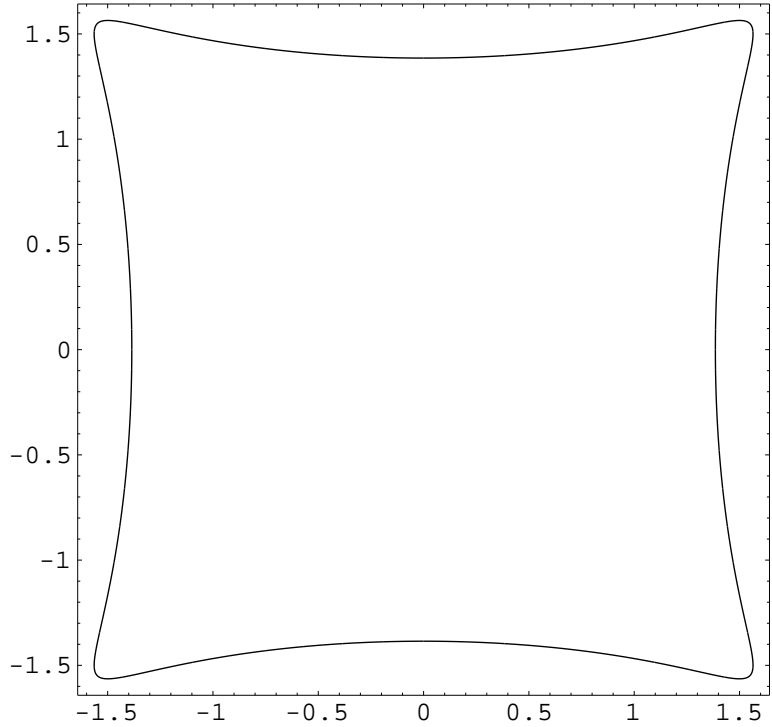}
    \caption{Spectrum of $u+u^*+i(v+v^*)$%, $1000\times 1000$ points
      }
  \end{center}
\end{figure}

\begin{Example}[$u+u^*+i(v+v^*)$]
  Let $u$ and $v$ be freely independent Haar unitaries (e.g., 
  $u=\lambda(g_1)$,   $v=\lambda(g_2)$, where $g_1$ and $g_2$ are the
  generators of the free group $\IF_2$ on two generators.
  The distribution of the selfadjoint operator $u+u^*$ is the arcsine
  distribution with Cauchy transform
  $G(\zeta)=\frac{1}{\zeta\sqrt{1-\frac{4}{\zeta^2}}}$.
  The operator $u+u^*+i(v+v^*)$ recently came up as a candidate for a
  counterexample to the invariant subspace conjecture in the von
  Neumann setting; however, its ``cousin'', the \emph{circular
  element} turned out to have plenty of invariant subspaces
  \cite{Dykema.Haagerup:2000:invariant}. Let us compute the border of
  its spectrum here.
  The moment generating function of $u+u^*$ is
  $$
  f(s) = \frac{1}{s}\,G(\frac{1}{s}) = \frac{1}{\sqrt{1-4s^2}}
  ;
  $$
  from this we easily obtain the moment generating function of $i(v+v^*)$
  $$
  g(t) = f(it) = \frac{1}{\sqrt{1+4t^2}}
  .
  $$
  Now the $L^2$-norm of the resolvent of a selfadjoint operator $a$ is
  \begin{align*}
    \norm{(1-sa)^{-1}}_2^2
    &= \frac{1}{\abs{s}^2}
       \norm{(\frac{1}{s}-a)^{-1}}_2^2\\
    &= - 
       \frac{1}{\abs{s}^2}
       \,
       \frac{G(\frac{1}{s})-G(\frac{1}{\bar s})}%
            {\frac{1}{s}-\frac{1}{\bar s}}\\
    &= \begin{cases}
         \frac{\bar s\, f(\bar s) - s\, f(s)}%
              {\bar s - s} & \Im s \ne 0\\
         f(s)+s\,f'(s)     & \Im s =0
       \end{cases}
  \end{align*}
  and we have the equations
  \begin{gather*}
  \lambda s f-f-g+1 = 0 \\
  s f-t g = 0 \\
  f^2 (1-4 s^2)-1 = 0 \\
  g^2 (1+4 t^2)-1 = 0 \\
  \end{gather*}
  A Gr\"obner basis with respect to degrevlex order on $(s,t,f,g)$ has
  $7$ elements and the quotient algebra $\IC[s,t,f,g]/I$ is four
  dimensional with basis $([g],[f],[t],[1])$. Elimination did not
  succeed in this example and we chose the eigenvalue approach of
  section~\ref{sub:ex:eigenvalueapproach}.  With the aid of the
  Gr\"obner basis the matrix of multiplication by $[s]$ is computed as
  \begin{equation*}
    M_s=
    \begin{bmatrix}
      \dfrac{4z^{2}-48}{z^{5}+24z^{3}+16z} & \dfrac{1}{z} & \dfrac{z^{2}+16}{z^{4}+24z^{2}+16} & \dfrac{z^{6}+16z^{4}-48z^{2}+128}{z^{7}+20z^{5}-80z^{3}-64z} \\
      \dfrac{z^{4}-24z^{2}-48}{z^{5}+24z^{3}+16z} & \dfrac{1}{z} & -\dfrac{2z^{2}-8}{z^{4}+24z^{2}+16} & -\dfrac{z^{6}-16z^{4}-48z^{2}-128}{z^{7}+20z^{5}-80z^{3}-64z} \\
      -\dfrac{z^{4}-24z^{2}+16}{z^{4}+24z^{2}+16} & 0 & \dfrac{3z^{3}+12z}{z^{4}+24z^{2}+16} & \dfrac{z^{6}-20z^{4}-80z^{2}+64}{z^{6}+20z^{4}-80z^{2}-64} \\
      \dfrac{z^{4}+16z^{2}+48}{z^{5}+24z^{3}+16z} & -\dfrac{1}{z} & -\dfrac{z^{2}+12}{z^{4}+24z^{2}+16} & \dfrac{8z^{2}+32}{z^{5}+24z^{3}+16z}
    \end{bmatrix}
  \end{equation*}
The result is shown in figure~\ref{fig:uuivv}.
\end{Example}

\begin{figure}[t]
  \label{fig:u2+u3+u4}
  \begin{center}
  \includegraphics[width=\textwidth, height=.4\textheight, keepaspectratio]{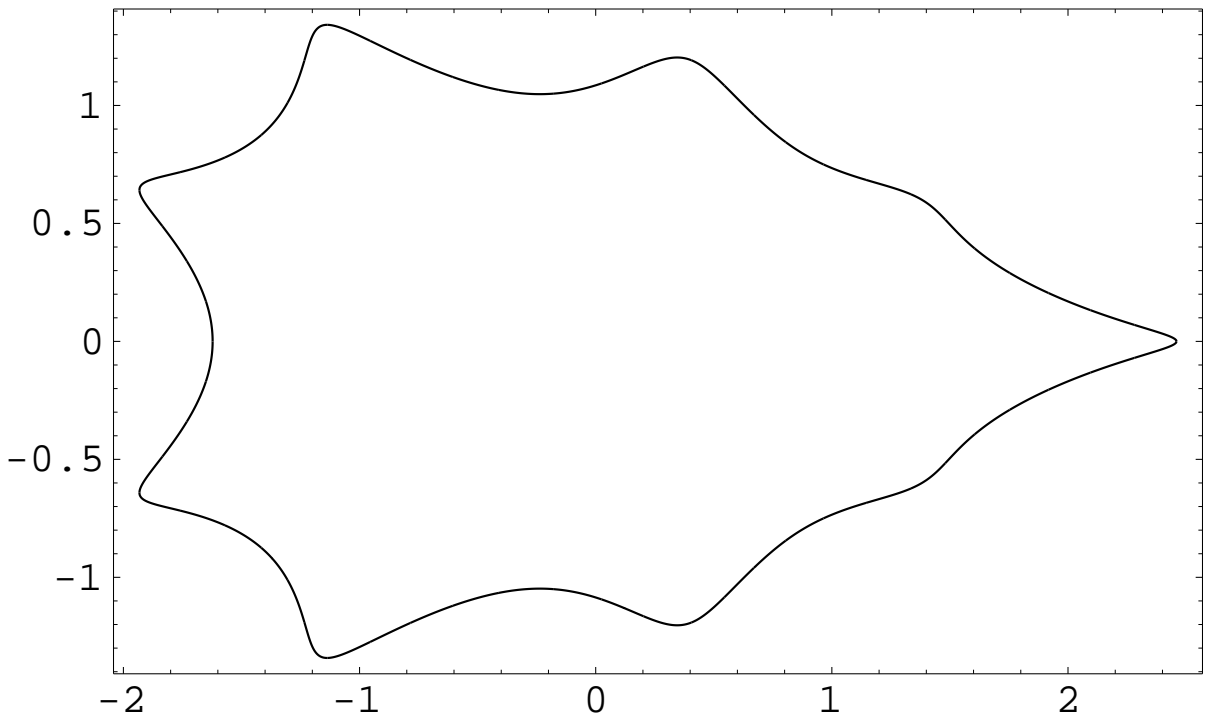}
    \caption{Spectrum of $u_2+u_3+u_4$%, $1000\times 600$ points
      }
  \end{center}
\end{figure}

\begin{Example}
  As an example with more than two summands, let us consider
  $u_2+u_3+u_4$.  Generally, the equations for a sum $\sum u_k$, where
  $u_k$ are free unitaries with $u_k^{n_k}=1$, are as follows.  If,
  for given $\lambda$, there is a solution $s_j$ for the system
  \begin{gather*}
    \lambda = \frac{1}{z} + \sum (\frac{1}{s_k} - \frac{1}{z})\\
    \frac{s_k}{1-s_k^{n_k}} = z \\
    \sum_k
     \frac{\abs{s_k}^2+\abs{s_k}^4+\dots+\abs{s_k}^{2n_k-2}}%
          {\abs{1-s_k^{n_k}}^2+\abs{s_k}^2+\abs{s_k}^4+\dots+\abs{s_k}^{2n_k-2}}
       < 1
  \end{gather*}
  then $\lambda$ is not in the spectrum.
  The example with $n_1=2$, $n_2=3$, $n_3=4$ is shown in
  fig.~\ref{fig:u2+u3+u4}.
\end{Example}

\emph{Acknowledgements.} We are grateful to U.~Haagerup and P.~Biane
for useful discussions. This work was done at Odense Universitet and
Universit\'e d'Orl\'eans under the EU-network ``Non-commutative
geometry'' ERB FMRX CT960073; the final part was done at the Institut
Henri Poincar\'e in Paris with a grant from the Ostrowski fundation.
We express our gratitude to all these institutions.

%BBBBBBBBBBBBBBBBBBBBBBBBBBBBBBBBBBBBBBBBBBBBBBBBBBBBBBBBBBBBBBBBBBBBBBBBB
\providecommand{\bysame}{\leavevmode\hbox to3em{\hrulefill}\thinspace}

\end{document}